\documentstyle[12pt]{article}
\titlepage
\newcommand{\ignore}[1]{}{}
\newcommand{\f}{\frac}

\renewcommand{\(}{\left(}
\renewcommand{\)}{\right)}
\renewcommand{\[}{\left[}

\newcommand{\no}{\nonumber}

\newcommand{\ep}{\epsilon}
\newcommand{\la}{\label}
\newcommand{\be}{\begin{eqnarray}}
\newcommand{\ee}{\end{eqnarray}}
\newcommand{\bestar}{\begin{eqnarray*}}
\newcommand{\eestar}{\end{eqnarray*}}

\newcommand{\al}{\alpha}

\newcommand{\IR}{I\!\! R}

\newcommand{\si}{\sigma}

\newcommand{\beq}{\begin{equation}}
\newcommand{\eeq}{\end{equation}}
\newcommand{\cD}{{\mathcal{D}}}
\newcommand{\VV}{\Vert\ \Vert}

\newtheorem{th}{Theorem}
\newtheorem{lem}{Lemma}
\newtheorem{cor}{Corollary}
\newtheorem{prop}{Proposition}

\begin{document}

\begin{center}
{\Large\bf On weighted approximations in $D[0, 1]$ with
applications
to self-normalized partial sum processes}\\
[0.8cm]
{\large Mikl\'os Cs\"org\H{o}} \\
Carleton University, Ottawa, Canada
\\[0.5cm]
{\large Barbara Szyszkowicz} \\
Carleton University, Ottawa, Canada\\
[0.5cm]
{\large Qiying Wang} \\
University of Sydney, Australia
\end{center}

\bigskip
\begin{center}
\textit {Dedicated to Istv\'an Berkes, S\'andor Cs\"org\H{o} and P\'eter Major in celebration of\\ their sixtieth years}
\end{center}

\bigskip
\centerline{{\bf ABSTRACT}}
\begin{center}
\begin{minipage}{120mm}

{Let $X, X_1, X_2,...$ be a sequence of non-degenerate i.i.d.
random variables with  mean zero. The best possible weighted
approximations are investigated in $D[0, 1]$ for the partial sum processes
$\{S_{[nt]}, 0\le t\le 1\}$, where $S_n=\sum_{j=1}^nX_j$, 
 under the assumption that $X$ belongs to the domain of attraction of the normal law.
 The conclusions then  are  used to establish similar results
  for the sequence of self-normalized partial sum processes
$\{S_{[nt]}/V_n, 0\le t\le 1\}$, where  $V_n^2=\sum_{j=1}^nX_j^2$.
$L_p$ approximations of self-normalized partial sum processes are
also  discussed.}

\end{minipage}
\end{center}



\medskip\noindent
{\bf Running Head:} Weighted approximations

\vskip2cm\noindent
--------------------------------------------

\noindent 
The research of M.~Cs\"org\H{o} and  B.~Szyszkowicz is supported by their  NSERC Canada Discovery Grants
at Carleton University, Ottawa, and Q.~Wang's research is supported in part by Australian Research Council at University of Sydney.

\newpage

\baselineskip=0.6 true cm

\section{Introduction  and main results
\la {sec1}}

Let $X, X_1, X_2,...$ be a sequence of non-degenerate i.i.d.
random variables with $EX=0$ and let $ S_n=\sum_{j=1}^nX_j$, $n\ge 1$,  denote their  partial sums. The
classical weak invariance principle in probability states that, on
an appropriate probability space for $\{S_n,\, n\geq 1\}$ and a standard Wiener process $\{W(t), \, 0\leq t<\infty\}$, as $n\to\infty$, we have
\be 
\sup_{0\le
t\le 1}\Big|\,\f {S_{[nt]}}{\sqrt n\sigma}-\f {W(nt)}{\sqrt n}
\,\Big|=o_P(1)\quad \mbox{if and only if}\quad Var
(X)=\sigma^2<\infty. 
\la {fini}
 \ee 
This invariance
principle in probability implies the following version of  Donsker's
classical Functional Central Limit Theorem on $(D,\cD,\Vert\ \Vert)$, where $\cD$ is the $\sigma$-field of subsets of $D=D[0,1]$ generated by its finite-dimensional subsets, and $\Vert \ \Vert$ stands for the sup-norm metric on $D$.

\medskip\noindent{\bf Corollary A} {\it With $\longrightarrow_d$ standing for convergence in distribution as $n\to\infty$, we have
$$
h(S_{[n\, \cdot]}/\sqrt n \delta) \longrightarrow_d h(W(\cdot))
$$
for all $h : D=D[0,1] \to \IR$ that are $(D,\cD)$ measurable and $\Vert \  \Vert$-continuous, or $\VV$-continuous except at points forming a set of Wiener measure zero on $(D,\cD)$, where $W(\cdot)$ stands for a standard Wiener process on the unit interval $[0,1]$.}

\medskip\noindent{\bf Remark A} \ For further reference the statement of Corollary A is summarized by writing, as $n\to\infty$,
$$
S_{[nt]}/\sqrt n \delta \Longrightarrow W(t) \hbox{ on } (D,\cD,\VV).
$$

For details and further thoughts on Donsker's theorem via (1) as in Corollary A, we refer to Breiman (1968, Sections 13.4--13.5), Cs\"org\H{o} and R\'ev\'esz (1981, Section 2.1), and Cs\"org\H{o} (2002, Sections 3.3--3.4).

In view of (\ref {fini}), it is natural to seek  conditions for
having the following weighted approximation on an appropriate probability space for $S_n$, $n\geq 1$, and $W(\cdot)$: 
\be 
\sup_{0< t\le
1}\Big|\,\f {S_{[nt]}}{\sqrt n\sigma}-\f {W(nt)}{\sqrt n}
\,\Big|\Big/q(t)&=&o_P(1), \quad n\to\infty, \la {in1} 
\ee 
where $q(t)$ is a non-negative function on $(0, 1]$ approaching  zero as
$t\downarrow 0$. In this direction, using the methods developed by
Cs\"org\H{o}, Cs\"org\H{o}, Horv\'ath and Mason [CsCsHM] (1986),
and subsequently by  Cs\"org\H{o} and Horv\'ath (1986), for
weighted empirical and quantile processes, based on the
Koml\'os, Major and Tusn\'ady [KMT] (1975, 1976)  and Major (1976)
strong approximations for partial sum processes, Cs\"org\H{o} and
Horv\'ath (1988) concluded that, on an appropriate probability
space, as $n\to \infty$, \be \sup_{1/n\le  t\le
1}n^{\mu}\,\Big|\,\f {{S_{[nt]}}}{\sqrt n\si}- \f {W(nt)}{\sqrt
n}\,\Big|\Big/t^{1/2-\mu}&=&O_P(1),\quad  \la {in2} \ee for any
$0\le \mu\le 1/2-1/r$, if we assume $E|X|^r<\infty$ for some
$r>2$.

Let $Q$  be the class of  positive functions $q(t)$ on $(0,1]$,
i.e., $\inf_{\delta\le t\le 1}q(t)>0$ for $0<\delta<1$,
 which are nondecreasing near zero, and let
$$
I(q,c)=\int_{0+}^1t^{-1}\exp\(-cq^2(t)/t\)dt, \quad 0<c<\infty.
$$
By virtue of the result (\ref {in2}), we have the following fact:
if $E|X|^r<\infty$ for some $r>2$ and  $q\in Q$ such that $I(q,
c)<\infty$ for all $c>0$, then the claim (\ref {in1}) holds true.
We note in passing that the class of weight functions $q\in Q$ such that $I(q, c)<\infty$ for all $c>0$
is the optimal class for having (\ref {in1}) (cf.~CsCsHM (1986) via Lemma 3 here).
For more details  along these lines, we refer to Szyszkowicz
(1991, 1992, 1996, 1997) where, by a different method,
 the statement of  (\ref {in1}) is established, again
for  the optimal class of weight functions,  assuming only the
existence of two moments.

In this paper, we  establish a result which is similar to (\ref
{in1}) only under the assumption that
 $X$ belongs to the domain of attraction of the normal law,
 written  $X\in$  DAN throughout.
 In fact, it is well-known that $X\in$  DAN with $EX=0$
 if and only if there exists a sequence of constants
 $d_n\uparrow \infty$ such that, as $n\to\infty$, $
S_n/{d_n}\, \longrightarrow_d N(0,1)$. It is natural to ask
whether a similar version of the weighted approximation (\ref
{in1}) could also hold true. Throughout the paper, we make use of
the notation $ l(x) = EX^2I_{(|X|\le x)}, $ $b=\inf\big\{x\ge 1:
l(x)>0\big\}$, \
$$
\eta_j=\inf \Big\{s: s\ge b+1,\frac  {l(s)}{s^2}\le \frac {1}
{j}\Big\},\quad j=1,2,...\,.$$ and $b_n^2=n\, l(\eta_n)$. We  have
the following theorem which provides an affirmative answer to this
paramount question (cf. Corollary \ref {cor1}).

\begin{th} \la {th1} Let $q\in Q$ and assume that  $X\in$  DAN and $EX=0$.
Then, on an appropriate probability space for $X,X_1, X_2,\cdots,$
a standard Wiener process $\{W(t), 0\le t<\infty\}$ can be
constructed in such a way that the following  statements hold true as $n\to\infty$:

\begin{itemize}
\item[\mbox{\rm (a)}]  If $I(q,c)<\infty$ for all $c>0$, then
\be \sup_{1/n\le t\le 1}\Big|\,S_{[nt]}/b_n-W(nt)/\sqrt
n\,\Big|\Big/q(t) &=&o_P(1). \la {adth11} \ee

\item[\mbox{\rm (b)}]  If $I(q,c)<\infty$ for some $c>0$, then
\be \sup_{1/n\le t\le 1}\Big|\,S_{[nt]}/b_n-W(nt)/\sqrt
n\,\Big|\Big/q(t) &=& O_P(1). \la {adth12} \ee

\item[{\rm (c)}]  If $I(q,c)<\infty$ for some $c>0$, then
there exists a sequence of constants $\tau_n\to 0$ such that \be
\sup_{\tau_n \le t\le 1}\Big|\,S_{[nt]}/b_n-W(nt)/\sqrt
n\,\Big|\Big/q(t) &=& o_P(1). \la {adth35} \ee
\end{itemize}
\end{th}

{\bf Remark 1.} If $EX^2<\infty$, then $\tau_n$ in (\ref {adth35})
can be changed to $1/n$. In this case, it is readily seen that
$n\si^2/b_n^2\to 1$, as $n\to\infty$. Hence, via  Theorem 1.1 of
Szyszkowicz (1997), a standard Wiener process $\{W(t), 0\le
t<\infty\}$ can be constructed such that
 \bestar
&&\sup_{1/n \le t\le 1}\Big|\,S_{[nt]}/b_n-W(nt)/\sqrt n\,\Big|\Big/q(t)\\
&&\qquad \le\ \f {\sqrt n\si}{b_n}
\sup_{1/n \le t\le 1}\Big|\,S_{[nt]}/{\sqrt {n}\si}-W(nt)/\sqrt n\,\Big|\Big/q(t)\\
&&\qquad \quad\qquad + \,\Big|\f {\sqrt {n}\si}{b_n}-1
\Big|\sup_{0<t\le 1}\Big|n^{-1/2}W(nt)\Big|\Big/q(t) \\
&&\qquad =\ o_P(1), \eestar where
  we make  use of
$ \sup_{0<t\le 1}\Big|n^{-1/2}W(nt)\Big|\Big/q(t)=O_P(1), $ which
follows from  Lemma \ref {alem3} below via having the equality in distribution
$$
\left\{n^{-1/2}W(nt)\Big/q(t), 0<t\le 1\right\} {=}_d \left\{W(t)\Big/q(t), 0<t\le 1\right\}
$$ 
for each
$n\geq 1$. However, it is impossible to replace $\tau_n$ in (\ref
{adth35}) by $1/n$ without further restrictions on $X$ and/or $q\in Q$.
 In Section \ref {sec2}, an example will be given [cf. Proof of (\ref {adth40})]
 to show that
(\ref {adth35}) holds in terms of  a standard Wiener process
$\{W(t), 0\le t<\infty\}$, but \be \sup_{1/n \le t\le
1}\Big|\,S_{[nt]}/b_n-W(nt)/\sqrt n\,\Big|\Big/q(t)\neq o_P(1),
\la {adth40} \ee with $q^2(t)=t\log\log (t^{-1})$, where $\log
x=\log (\max\{e, x\})$ here, as well as throughout.

\medskip

The following  corollaries are  consequences of Theorem \ref
{th1}, and are of independent interest.

\begin{cor} \la {cor1} Let $q\in Q$.
As $n\to \infty$, the following statements are equivalent:
\begin{itemize}
\item[{\rm(a)}] $X\in$ DAN and $EX=0$;

\item[{\rm  }] $S_{[nt]}/b_n \Rightarrow W(t)$ on $(D,{\cD},\Vert~/q\Vert)$
if and only if $I(q,c)<\infty$ for all $c>0$, where $\{W(t),\,
0\leq t\leq 1\}$ is a standard Wiener process and $\Vert ~/q\Vert$
is the weighted sup-norm metric for functions $x,y$ 
 in $D[0,1]$, defined  by
\beq 
||(x-y)/q|| = \sup_{0\leq t\leq 1} |(x(t)-y(t))/q(t)|
\label{eq3.7} 
\eeq 
whenever this is well defined, i.e., {\it when
$\limsup_{t\downarrow 0} |(x(t)-y(t))/q(t)|$ is finite};

\item[{\rm(c)}] On the probability space of {\rm Theorem \ref {th1}} for
$X,X_1,X_2,\ldots$, a standard Wiener process $\{W(t),\, 0\leq
t<\infty\}$ can be constructed in such a way that, as
$n\to\infty$, 
\beq 
\sup_{0<t\leq 1} \big| S_{[nt]}/b_n -
W(nt)/n^{1/2}\big|/q(t) = o_P(1) \label{eq3.10} 
\eeq 
if and only
if $I(q,c)<\infty$ for all $c>0$.
\end{itemize}
\end{cor}

\noindent Mutatis mutandis, (b) above is to be read {\it \`a la} Corollary A and Remark A.

\medskip
While Corollary \ref {cor1} is a direct consequence of Theorem
\ref {th1}, we now state convergence in distribution results for
sup-functionals of weighted normalized partial sums for the
optimal class of weight functions $q\in Q$ satisfying
$\lim\sup_{t\downarrow 0}|W(t)|/q(t)<\infty~ a.s.$ (see (b) of Lemma \ref
{alem3}).
 Consequently, the results that follow are not implied by Corollary \ref {cor1} above,
 and they cannot be obtained via classical methods of weak convergence either, for
 tightness in our weighted sup-norm is not guaranteed by (b) of Lemma \ref {alem3}.

\begin{cor} \la{cor3.4}
Let $\{W(t), 0\le t< \infty\}$ be a standard Wiener process, $X
\in$ DAN and $EX=0$.

\begin{itemize}

\item[{\rm (a)}] If  $q\in Q$, then there exists a sequence of
positive constants $\tau_n\to 0$ such that,
 as $n\to\infty$,
\beq b_n^{-1}\sup_{\tau_n\le t\le 1}|S_{[nt]}|/q(t) \ \to_d \
\sup_{0< t\le 1} |W(t)|/q(t) \la{eq3.12} \eeq if and only if
$I(q,c)< \infty$ for some $c> 0$.

\item[{\rm (b)}] If $q\in Q$ and  $q(t)$ is nondecreasing on $(0,1]$, then, as
$n\to\infty$,
\beq b_n^{-1}\sup_{0<  t\le 1}|S_{[nt]}|/q(t)\ \to_d\ \sup_{0<
t\le 1} |W(t)|/q(t) \la{eq3.13} \eeq if and only if $I(q,c)<
\infty$ for some $c> 0$. Consequently, as $n\to\infty$,
 we have
$$
\quad b_n^{-1}\sup_{0<t\le 1}|S_{[nt]}|\big/\big(t\log\log
(t^{-1})\big)^{1/2} \ \to_d\ \sup_{0< t\le 1}
|W(t)|\big/\big(t\log\log (t^{-1})\big)^{1/2}.
$$
\end{itemize}
\end{cor}

\medskip
{\bf Remark 2.} Corollary \ref {cor3.4} rhymes with Corollary 1.1
of  Szyszkowicz (1997), where $X$ is assumed to
 have two moments.
It remains an open problem whether $\sup_{\tau_n\le t\le 1}$ in
(\ref {eq3.12})
 can be improved to
$\sup_{0< t\le 1}$ without further restrictions on $X$ and/or
$q\in Q$.

\medskip
{\bf Remark 3.} The main results of Corollary \ref {cor1} were
announced without proofs
 in Cs\"org\H{o},
Szyszkowicz and Wang [CsSzW] (2004), where we reviewed  weighted
approximations and strong limit theorems for self-normalized
partial sum processes. It is interesting and also of interest to
note that the class of  the weight functions in Corollary \ref
{cor3.4} is bigger than that in Corollary \ref {cor1}. Such a
phenomenon was first noticed and proved for weighted empirical and
quantile processes by CsCsHM (1986) and then by Cs\"org\H{o} and
Horv\'ath (1988) for partial sums on assuming $E|X|^v<\infty$ for
some $v>2$. For more details along these lines, we refer to
Szyszkowicz
 (1991, 1996, 1997),
and to Cs\"org\H{o}, Norvai\v{s}a and Szyszkowicz [CsNSz](1999).

\medskip
We next consider  applications of Theorem \ref {th1} to the
so-called self-normalized partial sum processes defined by
$\{S_{[nt]}/V_n, 0\le t<\infty\}$, where
$V_n^2=\sum_{j=1}^nX_j^2$. Note that \be V_n^2/b_n^2 &\to_P& 1 \la
{impor3} \ee if $X\in$ DAN (see the result (18) in  CsSzW (2003),
for example). This, together with Theorem \ref {th1}, yields the
following Theorem \ref {th2}. We mention that Theorem \ref {th2}
 and its corollaries
 were announced without proofs in CsSzW (2004).
 We restate them here for convenient reference and  further use in the sequel.

\begin{th} \la {th2}
Assume that  $X\in$  DAN and $EX=0$. Then, on the  probability
space of {\rm Theorem \ref {th1}} for $X,X_1, X_2,\cdots,$ the
there constructed  standard Wiener process $\{W(t), 0\le
t<\infty\}$ is   such that that the following  statements hold
true.

\begin{itemize}
\item[{\rm (a)}] Let $q\in Q$. Then, as $n\to\infty$,
\be
&&\sup_{1/n\le t\le 1}\Big|\,S_{[nt]}/V_n-W(nt)/\sqrt n\,\Big|\Big/q(t)\no\\
&&\qquad \quad =\ \left\{\begin{array}{ll}
o_P(1) & \mbox{if $I(q,c)<\infty$ for all $c>0$},\\
O_P(1) & \mbox{if $I(q,c)<\infty$ for some $c>0$} .
\end{array}
\right. \la {adth1a} \ee

\item[{\rm (b)}]  Let $q\in Q$ and $I(q,c)<\infty$ for some $c>0$. Then
there exists a sequence of constants $\tau_n\to 0$ such that as
$n\to\infty$, \be \sup_{\tau_n \le t\le
1}\Big|\,S_{[nt]}/V_n-W(nt)/\sqrt n\,\Big|\Big/q(t)=o_P(1). \la
{adth35a} \ee
\end{itemize}
\end{th}

\begin{cor} \la {cor3} Let $q\in Q$.
As $n\to \infty$, the following statements are equivalent:
\begin{itemize}
\item[{\rm(a)}] $X\in$ DAN and $EX=0$;

\item[{\rm (b)}] $S_{[nt]}/V_n \Rightarrow_d W(t)$ on $(D,\cD,\Vert~/q\Vert)$
if and only if $I(q,c)<\infty$ for all $c>0$, where $\Vert
~/q\Vert$ is defined as in {\rm Corollary \ref {cor1}} and
 $\{W(t),\, 0\leq t\leq 1\}$ is a standard Wiener process;

\item[{\rm(c)}] On the probability space of {\rm Theorem \ref {th1}} for
$X,X_1,X_2,\ldots$, a standard Wiener process $\{W(t),\, 0\leq
t<\infty\}$ can be constructed in such a way that, as
$n\to\infty$, \beq \sup_{0<t\leq 1} \big| S_{[nt]}/V_n -
W(nt)/n^{1/2}\big|/q(t) = o_P(1) \label{eq3.10a} \eeq if and only
if $I(q,c)<\infty$ for all $c>0$.\end{itemize}
\end{cor}

\begin{cor} \la{cor4} The  conclusions of {\rm Corollary \ref{cor3.4}} continue to
 hold true if we replace
the normalizing constants $b_n$ by $V_n$.
\end{cor}

{\bf Remark 4.} With $q=1$, Corollary \ref {cor3} reduces to  Theorem 1 of CsSzW
(2003),  an  in probability version of
 Donsker's theorem for self-normalized partial sum processes.
 Corollaries \ref {cor3} and \ref {cor4} provide  basic tools
for investigating the limit behavior  of  statistics that arise in
studying the problem of a change in the mean in the domain of
attraction of the normal law. For  details we  refer to Section 5
in CsSzW (2004).

\medskip
For the sake of Studentized versions of Theorem \ref {th2},
Corollaries \ref {cor3} and \ref {cor4}, consider now the sequence
$T_{n,t}(.)$ of Student processes in $t\in [0,1]$ on $D[0,1]$,
defined as 
\bestar 
\{T_{n,t}(X), 0\le t\le 1\} &:=& \Bigg\{\f
{(1/\sqrt n)\, \sum_{i=1}^{[nt]}X_i}{ \sqrt {[1/(n-1)]\,
\sum_{i=1}^n(X_i-\bar X_n)^2}}, \ 0\le t\le 1\Bigg\} \no\\
&=&\Bigg\{\f {\, \sum_{i=1}^{[nt]}X_i/V_n}{ \sqrt
{\,[n-(S_n/V_n)^2]/(n-1)}}, \ 0\le t\le 1\Bigg\}. 
\eestar 
Clearly, $T_{n,1}(X)$ is the familiar form of the classical Student ratio.
When $X=_d N(\mu,\si^2)$, then $T_{n, 1}(X-\mu)$ is his famous
$t$-random variable with $n-1$ degrees of freedom [cf. ``Student"
(1908)]. Clearly, if $T_{n,1}(X)$ or $S_n/V_n$ has an asymptotic
distribution, then so does the other, and it is well known that
they coincide (cf., e.g., Proposition 1 in Griffin (2002)).
Hence, usually, and without loss of generality, the limiting
distribution of $S_n/V_n$ is studied in lieu of that of
$T_{n,1}(X)$.  In fact via Theorem 2 and Corollaries 3 and 4 we
follow a functional version of this route.  Backtracking a bit,
Logan, Mallows, Rice and Shepp (1973) conjectured that ``$S_n/V_n$
is asymptotically normal if (and perhaps only if) $X$ is in the
domain of attraction of the normal law". It is the ``only if" part
that has remained open until 1997 for the general case of not
necessarily symmetric random variables, when Gine, G\"otze and
Mason (1997) proved the following result.

{\bf Theorem A.} {\it The following two statements are equivalent:

\noindent
{\rm(i)}  ~$X$ is in the domain of attraction of the normal law;

\noindent
{\rm(ii)} There exists a finite constant $\mu$ such that, as $n\to\infty$,
$T_{n,1}(X-\mu)\to_d N(0,1)$. 

\smallskip
Moreover, if either {\rm(i)} or {\rm(ii)} holds, then $\mu=EX$.}

\medskip
Furthermore, Chistyakov and G\"otze (2004)  confirmed a second
conjecture of Logan et al. (1973) that the Student $t$-statistic
has a non-trivial limiting distribution if and only if $X$ is in
the domain of attraction of a stable law with some exponent
$\al\in (0,2]$. Our next corollary constitutes a Studentized
version of the self-normalized results of Theorem \ref {th2}, and
Corollaries \ref {cor3} and \ref {cor4}. Thus it amounts to
various optimal functional extensions of Theorem A in $(D,\cD,
\|/q\|)$.

\begin{cor} \la {cor5} The results {\rm (a), (b)} of {\rm Theorem \ref {th2}}
continue to
hold true when $EX=\mu$, with $T_{n, t}(X-\mu)$ replacing
$S_{[nt]}/V_n$ in each of the two statements. Mutatis mutandis,
the same is true as well in case of {\rm Corollaries \ref {cor3}
{\it and} \ref {cor4}}.
\end{cor}

\medskip
This paper is organized as follows. In Section 2 we provide proofs
of the  main results. The proof of Proposition \ref {pro1}, which
is used in the proofs of main results, will be given in  Section
\ref {sec3}. Finally,  in Section \ref{sec5} we will discuss
  weighted approximations of self-normalized  partial sum processes
under the optimal class of weight functions in $L_p$. Throughout
the paper, $A, A_1, A_2, ...$ denote constants which may be
different at each appearance.

\medskip
In concluding this section, we mention some  of the previous
results on weighted approximation. The study of weighted
approximations originates
 from the earlier works of Anderson and Darling (1952),  R\'enyi
(1953), Chibisov (1964), Pyke and Shorack (1968) and O'Reilly
(1974), who investigated the asymptotic behaviour of weighted
empirical and quantile processes. The works of these authors  were
later extented by CsCsHM (1986), Szyszkowicz (1991, 1992, 1996,
1997) and CsNSz (1999). In CsCsHM (1986), the authors established approximations
of empirical and quantile processes by sequences of Brownian
bridges in weighted supremum metrics for the optimal class of
weight functions. Szyszkowicz (1996, 1997)
 derived similarly optimal  weighted approximations of standardized partial sum
processes
 in $D[0,\infty)$ under finite second moment conditions
 (cf. also Szyszkowicz (1991) concerning related results in $D[0,1]$).
The paper CsNSz (1999) extends the results of Szyszkowicz (1991,
1992, 1996, 1997) to arbitrary positive weight functions $q$ on
$(0,\infty]$ for which $\lim_{t\downarrow 0, t\uparrow
0}|W(t)|/q(t)=0\, a.s.$, or $\limsup_{t\downarrow 0, t\uparrow
0}|W(t)|/q(t)<\infty\, a.s.$ For further results on weighted
approximations for empirical, quantile  and standardized partial
sum processes with applications to change-point analysis, we refer
to  two books of Cs\"org\H{o} and Horv\'ath (1993, 1997).

\section{Proofs of the main results
\la {sec2}}

We first list several lemmas of independent interest that will be
used in the proofs of the  main results of Section \ref {sec1}.
The first result is due to Sakhanenko (1980, 1984, 1985).

\begin{lem} \la {lem2} Let $ X_1, X_2,...$ be independent random variables
with $EX_j=0$  and $\sigma_j^2=EX_j^2<\infty $ for each $j\ge 1$.
Then we can redefine $\{X_j, j\ge 1\}$ on a richer  probability
space together with
 a sequence of independent $N(0,1)$ random variables, $Y_j,\, j\ge 1$,  such that
 for every $p>2$ and $x>0$,
\bestar P\Big\{\max_{i\le
n}\Big|\sum_{j=1}^iX_j-\sum_{j=1}^i\sigma_jY_j\Big|\ge x\Big\} \le
(Ap)^px^{-p}\sum_{j=1}^nE|X_j|^p, \eestar where $A$ is an absolute
positive constant.
\end{lem}

\medskip
The next two lemmas are due to CsCsHM (1986) [cf. Lemma A.5.1 and
Theorem A.5.1  respectively in Cs\"org\H{o} and Horv\'ath (1997)].
Proofs of Lemmas \ref {lem3} and \ref {alem3} can also be found in
Section 4.1 of Cs\"org\H{o} and Horv\'ath (1993). For related
further results along these lines we refer to Cs\"org\H{o}, Shao
and Szyszkowicz (1991) and, for some historical references, to
Section 2 of CsNSz (1999) (see the statements (2.5), (2.9) and
their discussion in there).

\begin{lem} \la {lem3} Let $q\in Q$.  If $I(q,c)<\infty$ for some $c>0$,
 then
$ \lim_{t\downarrow 0}\,t^{1/2}/q(t)=0. $

\end{lem}

\begin{lem} \la {alem3} Let  $\{W(t), 0\le t<\infty\}$
be a standard Wiener process and  $q \in Q$. Then,
\begin{itemize}
\item[{\rm (a)}] $I(q,c)<\infty$ for all $c>0$ if and only if
$ \limsup_{t\downarrow 0}\,|W(t)|/q(t)=0,\ a.s.$

\item[{\rm (b)}] $I(q,c)<\infty$ for some $c>0$ if and only if
$ \limsup_{t\downarrow 0}\,|W(t)|/q(t)<\infty,\
 a.s. $
\end{itemize}
\end{lem}

\medskip
We are now ready to prove the main results of Section \ref {sec1}.

\textbf{Proof of Theorem \ref {th1}.} In addition to the notation
in Section \ref {sec1}, write \bestar
  \quad X_j^*=X_jI_{(|X_j|\le \eta_j)}\qquad \mbox{and}\qquad
  S_n^*=\sum_{j=1}^nX_j^*.
\eestar By  Lemma \ref {lem2}, we can redefine $\{X_j, j\ge 1\}$
on a richer  probability space together with
 a sequence of independent $N(0,1)$ random variables, $Y_j,\, j\ge 1$,
 such that for any $x>0$ and  any constant sequence $c_j, j\ge 1$,
\be P\Big\{\max_{i\le n}\Big|\sum_{j=1}^ic_j\(X_j^*-EX_j^*\)-
\sum_{j=1}^ic_j\sigma_j^*Y_j\Big|\ge x \Big\} \ \le\
Ax^{-3}\sum_{j=1}^n|c_j|^{3}E|X|^3I_{(|X|\le \eta_j)}, \la {main1}
\ee where $\sigma_j^{*2}=Var (X_j^*)$. Let  $\{W(t), 0\le
t<\infty\}$ be a standard Wiener process such that
$$
W(n)=\sum_{j=1}^nY_j, \qquad n=1,2,3,...
$$
The results (\ref {adth11})-(\ref {adth35}) will be shown to hold
true for the this way constructed Wiener process, and then Theorem
\ref {th1} follows accordingly. To prove (\ref {adth11})-(\ref
{adth35}), we need the following proposition.

\begin{prop} \la {pro1} We have, as $n\to\infty$,
\be \f 1n\sum_{j=1}^n\Big(\f
{\sigma_j^*}{l^{1/2}(\eta_n)}-1\Big)^2 &\to& 0,
\la {ad188} \\
b_n^{-1} \sup_{0\le t\le 1}
\Big|S_{[nt]}-\sum_{j=1}^{[nt]}\si_j^*Y_j\Big|&=&o_P(1),
\la {formu4a} \\
b_n^{-1} \sup_{0< t\le 1}
\Big|S_{[nt]}-\sum_{j=1}^{[nt]}\si_j^*Y_j\Big|\Big/t^{1/2}&=&O_P(1).
\la {formu4} \ee
\end{prop}

The proof of Proposition \ref {pro1} will be given in Section 4.
By virtue of  Proposition \ref {pro1}, we have that if $I(q,
c)<\infty$ for some $c>0$, then \be I_n&:=&b_n^{-1} \sup_{1/n\le
t\le 1} \Big|S_{[nt]}-\sum_{j=1}^{[nt]}\si_j^*Y_j\Big|\Big/q(t)\
=\ o_P(1). \la {a17} \ee Indeed it follows from (\ref {formu4a})
and (\ref {formu4}) that, for any $0<\delta<1$, \bestar I_n &\le &
\sup_{1/n\le t\le \delta}t^{1/2}/q(t)\,\,b_n^{-1} \sup_{0< t\le 1}
\Big|S_{[nt]}-\sum_{j=1}^{[nt]}\si_j^*Y_j\Big|\Big/t^{1/2} \no\\
&&\qquad \quad\qquad   +\  \sup_{\delta<t\le
1}q^{-1}(t)\,\,b_n^{-1}\sup_{\delta<t\le 1}
 \Big|S_{[nt]}-\sum_{j=1}^{[nt]}\si_j^*Y_j\Big|\no\\
 &=& O_P(1)\sup_{1/n\le t\le \delta}t^{1/2}/q(t) +o_P(1) \sup_{\delta<t\le 1}q^{-1}(t).
\eestar The claim (\ref {a17}) now follows from Lemma \ref {lem3}
and the fact that $\inf_{\delta\le t\le 1}q(t)>0$ for any
$\delta>0$.

We now proceed to prove (\ref {adth11})-(\ref {adth35}). Consider
(\ref {adth11}) and (\ref {adth12}) first.   We have \be
 \sup_{1/n\le t\le 1}\Big|S_{[nt]}/b_n-W(nt)/\sqrt n\Big|\Big/q(t) &\le& I_n
+I_1(n)+I_2(n),
\la {formu1} \ee where $I_1(n)=n^{-1/2}\,\sup_{1/n\le t\le
1}\Big|W([nt])-W(nt)\Big|\Big/q(t)$ and
$$
I_2(n)=\sup_{1/n\le t\le
1}\Big|b_n^{-1}\sum_{j=1}^{[nt]}\sigma_j^*Y_j-n^{-1/2}
W([nt])\Big|\Big/q(t).
$$
Note that $\left\{\f {W(nt)-W([nt])}{\sqrt {nt}}, 0<t\le
1\right\}{=}_d \left\{W\big(\f
{nt-[nt]}{{nt}}\big), 0<t\le 1\right\}$ and $\f
{|nt-[nt]|}{{nt}}\le 1$ for $1/n\le t\le 1$. Similarly to the
proof of (\ref {a17}), it follows that 
\be 
I_1(n) &\le&
\sup_{1/n\le t\le \delta}t^{1/2}/q(t)\,
\sup_{1/n\le t\le \delta}\Big|W([nt])-W(nt)\Big|\Big/\sqrt {nt}\no\\
&&\qquad +
 n^{-1/2}\,\sup_{\delta\le t\le 1}q^{-1}(t)\, \sup_{\delta\le t\le 1}
 \Big|W([nt])-W(nt)\Big|\no\\
 &=&o_P(1), \la {a18}
\ee 
whenever $I(q, c)<\infty$ for some $c>0$. As to $I_2(n)$, we
have \be I_2(n) &\le & \sup_{1/n\le t\le
\delta}\Big|b_n^{-1}\sum_{j=1}^{[nt]}\sigma_j^*Y_j-n^{-1/2}
W([nt])\Big|\Big/q(t)  + \,II(n)\,\sup_{\delta\le t\le
1}q^{-1}(t), \la{a19} \ee for any $\delta\in (0,1)$, where
$II(n)=\, \max_{1\le k\le n}
\Big|b_n^{-1}\sum_{j=1}^{k}\sigma_j^*Y_j-
n^{-1/2}\sum_{j=1}^{k}Y_j\Big|$.  Furthermore,  $II(n)=o_P(1)$,
 since it follows from (\ref {ad188}) that, as $n\to\infty$,
\bestar 
E[II(n)]^2 & \le&
 \f 1n\, E\max_{1\le k\le n}\Big|\sum_{j=1}^{k}\Big(\f
{\sigma_j^*}{l^{1/2}(\eta_n)}-1\Big)Y_j\Big|^2 \le
\f 1n\sum_{j=1}^n\Big(\f {\sigma_j^*}{l^{1/2}(\eta_n)}-1\Big)^2 \
\to\ 0. 
\eestar

In view of  (\ref {a17})-(\ref {a19}), the results (\ref {adth11})
and (\ref {adth12}) will follow if we prove \be I_2^{(1)}(n,
\delta) &:=& \sup_{1/n\le t\le
\delta}\Big|b_n^{-1}\sum_{j=1}^{[nt]}\sigma_j^*Y_j-n^{-1/2}
W([nt])\Big|\Big/q(t)  \no\\
&=& \left\{\begin{array}{ll}
o_P(1) & \mbox{if $I(q,c)<\infty$ for all $c>0$},\\
O_P(1) & \mbox{if $I(q,c)<\infty$ for some $c>0$} ,
\end{array}
\right. \la {a12} \ee by letting $n\to\infty$ first and then
$\delta\downarrow 0.$ In order to prove (\ref {a12}), let
$\delta>0$ be small enough,
 so that $q(t)$ is already nondecreasing on $(0,\delta)$
and let $n$ be such that $1/n<\delta$. Since $\sigma_j^*\le
\Big(EX^2I_{(|X|\le \eta_j)}\Big)^{1/2}\le l^{1/2}(\eta_n), 1\le
j\le n,$ we have $\f 1n\sum_{j=1}^{[nt]}\Big(\f
{\sigma_j^*}{l^{1/2}(\eta_n)}-1\Big)^2\le t$, for $ t\in [0,1]$.
On the other hand, for each $n\ge 1$,
$$
\left\{n^{-1/2}\sum_{j=1}^{[nt]} \Big(\f
{\sigma_j^*}{l^{1/2}(\eta_n)}-1\Big)Y_j, 0<t\le 1\right\}\stackrel
{\rm {\cal D}}{=} \left\{W\Big(\f 1n\sum_{j=1}^{[nt]}\Big(\f
{\sigma_j^*}{l^{1/2}(\eta_n)}-1\Big)^2\Big), 0<t\le 1\right\}.
$$
Now it is readily seen  that \be I_2^{(1)}(n, \delta) &\stackrel
{\rm {\cal D}}{=} & \sup_{1/n\le t\le \delta}\Big|W\Big(\f
1n\sum_{j=1}^{[nt]}
\Big(\f {\sigma_j^*}{l^{1/2}(\eta_n)}-1\Big)^2\Big)\Big|\Big/q(t) \no\\
&\le & \sup_{1/n\le t\le \delta}\Big|W\Big(\f 1n\sum_{j=1}^{[nt]}
\Big(\f {\sigma_j^*}{l^{1/2}(\eta_n)}-1\Big)^2\Big)\Big|\Big/
q\Big(\f 1n\sum_{j=1}^{[nt]}\Big(\f {\sigma_j^*}{l^{1/2}(\eta_n)}-1\Big)^2\Big)\no\\
&\le& \sup_{0< t\le \delta}\Big|W(t)|/q(t). \la {adfor3} \ee This,
together with  Lemma  \ref {alem3},  implies (\ref {a12}). The
proofs of (\ref {adth11}) and (\ref {adth12}) is now complete.

We next prove (\ref {adth35}). Similarly to the proofs of (\ref
{adth11}) and (\ref {adth12}), it suffices to show that there
exists a sequence of positive constants $\tau_n\to 0$ such that if
$I(q,c)<\infty$ for some $c>0$, then \be \sup_{\tau_n\le t\le
\delta}\Big|b_n^{-1}\sum_{j=1}^{[nt]}\sigma_j^*Y_j-n^{-1/2}
W([nt])\Big|\Big/q(t)=o_P(1), \la {a15} \ee when $n\to\infty$
first and then $\delta\downarrow 0$. In fact,  in view of (\ref
{ad188}), there exists a sequence of positive constants
$\kappa_n\to\infty$ so that
$$
\f {\kappa_n}n\,\sum_{j=1}^n\Big(\f
{\sigma_j^*}{l^{1/2}(\eta_n)}-1\Big)^2 \to 0, \quad  \mbox{as
$n\to\infty$}.
$$
This implies that for any $\ep>0$,
$$
\f {1}{nt}\,\sum_{j=1}^{[nt]}\Big(\f
{\sigma_j^*}{l^{1/2}(\eta_n)}-1\Big)^2 \le \f
{\kappa_n}n\,\sum_{j=1}^n\Big(\f
{\sigma_j^*}{l^{1/2}(\eta_n)}-1\Big)^2 \le \ep^2,
$$
whenever $ 1/\kappa_n\le t\le 1$ and $n$ is large enough. Let
$\tau_n=1/\kappa_n$ and $\delta$ be small enough so that $q(t)$ is
nondecreasing on $(0,\delta)$. Then $\tau_n\to 0$ and, similarly
to the proof of (\ref {adfor3}), we have that if $I(q,c)<\infty$
for some $c>0$, then 
\be 
&&\sup_{\tau_n\le t\le
\delta}\Big|n^{-1/2}\sum_{j=1}^{[nt]} \Big(\f
{\sigma_j^*}{l^{1/2}(\eta_n)}-1\Big)Y_j\Big|\Big/q(t)
\no\\
&{=}_d & \sup_{\tau_n\le t\le
\delta}\Big|W\Big(\f 1n\sum_{j=1}^{[nt]}
\Big(\f {\sigma_j^*}{l^{1/2}(\eta_n)}-1\Big)^2\Big)\Big|\Big/q(t)\no\\
&{=}_d & \sup_{\tau_n\le t\le \delta} \f
{\Big|\ep\, W\Big(\f 1{n\ep^2}\sum_{j=1}^{[nt]} \Big(\f
{\sigma_j^*}{l^{1/2}(\eta_n)}-1\Big)^2\Big)\Big|} {q\Big(\f
1{n\ep^2}\sum_{j=1}^{[nt]}\Big(\f
{\sigma_j^*}{l^{1/2}(\eta_n)}-1\Big)^2\Big) }\, \f {q\Big(\f
1{n\ep^2}\sum_{j=1}^{[nt]}\Big(\f
{\sigma_j^*}{l^{1/2}(\eta_n)}-1\Big)^2\Big)
}{q(t)}\no\\
&\le &\ep\, \sup_{0< t\le \delta}\Big|W(t)|/q(t)=\ep\, O_P(1). \la
{adfor9} \ee This yields (\ref {a15}) and hence the proof of (\ref
{adth35}). The proof of Theorem \ref {th1} is now complete.

\medskip

\textbf{ Proof of Corollary \ref {cor1}.} We only show part (c)
following from part (a). Clearly, (c) implies (b), and that, in
turn, (a). In fact, under the conditions that $X\in$ DAN and
$EX=0$, it follows from Theorem~\ref {th1} and  Lemma \ref {alem3}
that \bestar
&&\sup_{0< t\le 1}\Big|S_{[nt]}/b_n-W(nt)/\sqrt n\Big|\Big/q(t)\\
&&\qquad \le \ \sup_{0<t<1/n}\Big|n^{-1/2}W(nt)\Big|\Big/q(t)+
\sup_{1/n\le t\le 1}\Big|S_{[nt]}/b_n-W(nt)/\sqrt n\Big|\Big/q(t)\\
&&\qquad = o_P(1), \eestar if $I(q,c)<\infty$ for any $c>0$. This
shows sufficiency of part $(\mbox{c})$. Noting that $ \sup_{0<
t\le 1}\Big|S_{[nt]}/V_n-W(nt)/\sqrt n\Big|\Big/q(t)\ge
\sup_{0<t<1/n}\Big|n^{1/2}W(nt)\Big|\Big/q(t), $ the necessity of
part $(\mbox{c})$ follow as in
 Szyszkowicz (1996) (pages 331-332). Hence we omit these details.
The proof of Corollary \ref {cor1} is now complete.

\medskip
\textbf{Proof of Corollary \ref {cor3.4}.} If one of (\ref
{eq3.12}) and (\ref {eq3.13}) holds, then \be \sup_{0<t\le 1}\Big|
W(t)/q(t)\Big|<\infty \quad a.s., \la{456} \ee for any standard
Wiener process. Consequently, Lemma \ref {alem3} implies that
$I(q,c)<\infty$ for some $c>0$.

Next we assume $I(q,c)<\infty$ for some $c>0$. By Lemma \ref
{alem3} again, we get (\ref {456}). Then,  clearly as
$n\to\infty$, we have that 
\bestar 
\sup_{\tau_n\le t\le
1}|W(t)|/q(t)\rightarrow_d \sup_{0< t\le 1} |W(t)|/q(t), 
\eestar
which, together with  part $(\mbox{c})$ of Theorem \ref {th1},
implies  (\ref {eq3.12}).

We next prove (\ref {eq3.13}). Recalling  (\ref {a17}), we have
\be b_n^{-1}\sup_{0< t\le 1}\Big|S_{[nt]}\Big|\Big/q(t) &\le&
\sup_{1/n\le t\le
1}\Big|b_n^{-1}\sum_{j=1}^{[nt]}\sigma_j^*Y_j\Big|\Big/q(t)
+o_P(1). \la {reformu1} \ee Since $\sigma_j^{*2}\le EX^2I_{(|X|\le
\eta_j)}=l(\eta_j)$, it can be easily seen that $
 b_n^{-2}\sum_{j=1}^{[nt]}
\sigma_j^{*2}\le t$ for $t\in (0,1]$. This, together with the fact
that
 $q(t)$ is nondecreasing  on $(0,1]$, yields that
\be 
\sup_{1/n\le t\le
1}\Big|b_n^{-1}\sum_{j=1}^{[nt]}\sigma_j^*Y_j\Big|\Big/q(t)
&{=}_d & \sup_{1/n\le t\le
1}\Big|W\Big(b_n^{-2}\sum_{j=1}^{[nt]}
\sigma_j^{*2}\Big)\Big|\Big/q(t)\no\\
&\le & \sup_{1/n\le t\le 1}\Big|W\Big(b_n^{-2}\sum_{j=1}^{[nt]}
\sigma_j^{*2}\Big)\Big|\Big/q\Big(b_n^{-2}\sum_{j=1}^{[nt]}
\sigma_j^{*2}\Big)\no\\
&\le& \sup_{0< t\le 1}\Big|W(t)|/q(t). \la {addtion} \ee It
follows from  (\ref {reformu1}) and (\ref {addtion})  that
 for any $x\ge 0$
$$
\lim_{n\to\infty}P\(b_n^{-1}\sup_{0< t\le
1}\Big|S_{[nt]}\Big|\Big/q(t) \le x\) \ge P\(\sup_{0< t\le
1}\Big|W(t)|/q(t)\le x\).
$$
On the other hand, (\ref {eq3.12}) and the fact that $\tau_n>0$
obviously imply that for any $x\ge 0$ \bestar
\lim_{n\to\infty}P\(b_n^{-1}\sup_{0< t\le
1}\Big|S_{[nt]}\Big|\Big/q(t) \le x\) &\le&
\lim_{n\to\infty}P\(b_n^{-1}\sup_{\tau_n\le t\le
1}\Big|S_{[nt]}\Big|\Big/q(t)
\le x\) \\
&=&
 P\(\sup_{0< t\le 1}\Big|W(t)|/q(t)\le x\).
\eestar Therefore,  we get the desired (\ref {eq3.13}). This also
completes the proof of Corollary \ref {cor3.4}.

\medskip
\textbf{Proof of (\ref {adth40})}.  Let us consider i.i.d.
symmetric random variables $X, X_1, X_2,...$ satisfying
$EX^2I_{(|X|\le x)}=0$, for $x\le 1$, and \be l(x)=EX^2I_{(|X|\le
x)} \sim \exp\((\log x)^{\alpha}\), \quad \mbox{for $x>1$}, \la
{tue} \ee where $0<\alpha<1$. Then $l(x)$ is a slowly varying
function at $\infty$. Hence,
 $EX=0$ and $X$ is in the domain of attraction of the normal
law. Along the proof,  and using the notations, of Theorem~\ref{th1}, a
standard Wiener process $\{W(t), 0\le t<\infty\}$ can be
constructed such that (\ref {adth35}) holds and \bestar
&&\sup_{1/n\le t\le 1}\Big|S_{[nt]}/b_n-W(nt)/\sqrt n\Big|\Big/q(t)\no\\
&&\qquad \ge\ \sup_{1/n\le t\le
1}\Big|b_n^{-1}\sum_{j=1}^{[nt]}\sigma_j^*Y_j-n^{-1/2}W(nt)\Big|\Big/q(t)
-I(n)\no\\
&&\qquad \ge\ \sup_{1/n\le t\le 1/\sqrt n}
\Big|b_n^{-1}\sum_{j=1}^{[nt]}\sigma_j^*Y_j-n^{-1/2}W(nt)\Big|\Big/q(t)-I(n)
\no\\
&&\qquad \ge\ \sup_{1/n\le t\le 1/\sqrt n}
\Big|n^{-1/2}W(nt)\Big|\Big/q(t)-\sup_{1/n\le t\le 1/\sqrt n}
\Big|b_n^{-1}\sum_{j=1}^{[nt]}\sigma_j^*Y_j\Big|\Big/q(t) -I(n),
\eestar where $I(n)$ is defined as in (\ref {a17}). Noting that
$q^2(t)=t\log\log (t^{-1})$ satisfies $I(q, 2)<\infty$,
$I(n)=o_P(1)$. So, to show (\ref {adth40}), it suffices to show
that \be \sup_{1/n\le t\le 1/\sqrt n}
\Big|n^{-1/2}W(nt)\Big|\Big/q(t){=}_d
\sup_{1/n\le t\le 1/\sqrt n} \Big|W(t)\Big|\Big/q(t) \neq o_P(1)
\la {00} \ee with $q^2(t)=t\log\log (t^{-1})$, and if
$I(q,c)<\infty$ for some $c>0$, then \be \sup_{1/n\le t\le 1/\sqrt
n}
\Big|b_n^{-1}\sum_{j=1}^{[nt]}\sigma_j^*Y_j\Big|\Big/q(t)=o_P(1).
\la{01} \ee

We first prove (\ref {01}). For the $l(x)$ defined in (\ref
{tue}), it can be easily shown that, for $n$ large enough,
$\eta_n\ge \sqrt n$ and $\max_{1\le j\le \sqrt n}\eta_j\le \(\sqrt
n\)^{3/5}$. Hence,
$$
\max_{1\le j\le \sqrt n}l(\eta_j)/l(\eta_n)\le 2\exp\Big[
\(0.3^{\alpha}-0.5^{\alpha}\)(\log n)^{\alpha}\Big] \to 0,\quad
\mbox{as $n\to\infty$}.
$$
This, together with $\sigma_j^{*2}\le l(\eta_j)$, yields that, as
$n\to\infty$,
$$
\f 1{tb_n^2}\sum_{j=1}^{[nt]}\sigma_j^{*2}\le \f
1{nt}\sum_{j=1}^{[nt]}\f {l(\eta_j)}{l(\eta_n)}\to 0,
$$
whenever $1/n\le t\le \sqrt n$. Now, a method that is similar to making statements as in 
(\ref {adfor9}) shows that
$$
\sup_{1/n\le t\le 1/\sqrt n}
\Big|b_n^{-1}\sum_{j=1}^{[nt]}\sigma_j^*Y_j\Big|\Big/q(t)
{=}_d \sup_{1/n\le t\le 1/\sqrt n}
\Big|W\Big(b_n^{-2}\sum_{j=1}^{[nt]}\sigma_j^{*2}\Big)\Big|\Big/q(t)
=o_P(1),
$$
which implies (\ref {01}).

Next we prove (\ref {00}). In fact, if (\ref {00}) is not true,
i.e., if
$$
\sup_{1/n\le t\le 1/\sqrt n} \Big|W(t)\Big|\Big/q(t) = o_P(1),
$$
then there exists an $0<\ep\le 1/4$ such that
$$
P\(\sup_{1/n\le t\le 1/\sqrt n} \Big|W(t)\Big|\Big/q(t) \le
\ep^{1/2}/2\)\ge 1/2.
$$
Hence, following the proof of (4.1.14) in Cs\"org\H{o} and
Horv\'ath (1993, page 185), we have, for $q^2(t)=t\log\log
(t^{-1})$, as $n\to\infty$, \bestar 0& \leftarrow &P\(\sup_{1/n\le
t\le 1/\sqrt n} \Big|W(t)\Big|\Big/q(t)
> \ep^{1/2}/2\)\no\\
&& \ge\ \f 14 \int_{1/n}^{1/\sqrt n}\f 1 t\exp\(-2\ep q^2(t)/t\)dt\no\\
&&\ge\ \f 14 \int_{1/n}^{1/\sqrt n}\f 1 t\exp\(-2^{-1} \log\log(t^{-1})\)dt\no\\
&&\ge \ \f 12\Big(1-\f 1{\sqrt 2}\Big)\(\log n\)^{1/2}\to\infty.
\eestar This is a contradiction, which implies that (\ref {00})
holds. The proof of (\ref {adth40}) is now complete.

\medskip
\textbf{Proof of Theorem \ref {th2}.} It follows from (\ref
{adth12}) that if $I(q, c)<\infty$ for some $c>0$, then \bestar
b_n^{-1}\sup_{1/n\le t\le 1} \big|S_{[nt]}\big|\Big/q(t) &\le&
\sup_{1/n\le t\le 1}
n^{-1/2}\big|W(nt)\big|\Big/q(t) \no\\
&&+\sup_{1/n\le t\le 1} \Big|S_{[nt]}/b_n- W(nt)/\sqrt
n\Big|\Big/q(t) =O_P(1). \eestar This, together with (\ref
{impor3}),
 yields that
 \be
\sup_{1/n\le t\le 1} \Big|S_{[nt]}/V_n-S_{[nt]}/b_n\Big|\Big/q(t)
&\le& \Big|\f {b_n}{V_n}-1\Big|\,b_n^{-1}\sup_{1/n\le t\le 1}
\big|S_{[nt]}\big|\Big/q(t)=o_P(1), \la {a34}
 \ee
 whenever $I(q, c)<\infty$ for some $c>0$. Theorem \ref {th2}  now follows
immediately from
 (\ref {a34}) and Theorem \ref {th1}.

\medskip
\textbf{Proofs of Corollaries \ref {cor3}-\ref {cor5}}. In view of
Theorem \ref {th2}, the proofs of Corollaries \ref {cor3}-\ref
{cor5} are the same as in the proofs of Corollaries \ref {cor1}
and \ref {cor3.4} and hence the details are omitted.

\section{Proof of Proposition \ref {pro1}
\la {sec3}}

We only prove (\ref {formu4}).  The result (\ref {ad188}) can be
found in CsSzW (2003, page 1235) and the result (\ref
{formu4a}) follows from similar arguments. Let
$$
Z_j=X_j^*-EX_j^*-\sigma_j^*Y_j,\quad j=1,2,...,
$$
and write $\eta_0=0$. Note that $\eta_j^2\le (j+1)l(\eta_j)$ and
$l(\eta_n)=\sum_{k=1}^nEX^2I_{(\eta_{k-1}<|X|\le \eta_k)}.$ It
follows from the Shorack and Smythe inequality (cf. Shorack and
Weller, 1986, p. 844), (\ref {main1}) with $c_j=1/j^{1/2}$  that,
 for any $ C>0$,
\bestar &&P\(\sup_{1/n\le t\le 1}
\Big|\sum_{j=1}^{[nt]}Z_j\Big|\Big/t^{1/2}\ge C\, b_n\)\no\\
&&\qquad\le\  P\(
\max_{1\le k\le  n}\Big|(k/n)^{-1/2}\sum_{j=1}^{k}Z_j\Big|\ge C\, b_n\)\no\\
&&\qquad\le\  P\(
\max_{1\le k\le  n}\Big|\sum_{j=1}^{k}j^{-1/2}Z_j\Big|\ge 2^{-1}C\, b_n/\sqrt n\)\no\\
&&\qquad \le\ A\(\f {2\sqrt n}{C\, b_n}\)^3\sum_{j=1}^{ n}j^{-3/2}E|X|^3I_{(|X|\le
\eta_j)}\no\\
&&\qquad  \le\ AC^{-3}l^{-3/2}(\eta_n)\sum_{k=1}^{n}k^{-1/2}
E|X|^3I_{(\eta_{k-1}<|X|\le \eta_k)}\no\\
&&\qquad  \le\
AC^{-3}l^{-3/2}(\eta_n)\sum_{k=1}^{n}l^{1/2}(\eta_k)
E|X|^2I_{(\eta_{k-1}<|X|\le \eta_k)}\no\\
&&\qquad \le \ AC^{-3}. \la {a1} \eestar This yields \be
J_{1n}&:=&b_n^{-1}\,\sup_{1/n\le t\le 1}
\Big|\sum_{j=1}^{[nt]}Z_j\Big|\Big/t^{1/2}=O_P(1). \la {formu410}
\ee Similarly, by noting that
$l(\eta_n)=\sum_{k=1}^nEX^2I_{(\eta_{k-1}<|X|\le \eta_k)}$ and
$E|X|I_{(|X|>\eta_n)}=o(\eta_n^{-1}l(\eta_n))=o(b_n/n)$ [see Lemma
1 of CsSzW (2003) for example], it follows from (7) in CsSzW
(2003) with $a_j=EX^2I_{(\eta_{j-1}<|X|\le \eta_j)}$ that \bestar
&&P\(\sup_{1/n\le t\le 1}
\Big|\sum_{j=1}^{[nt]}(X_j-X_j^*+EX_j^*)\Big|\Big/t^{1/2}\ge C\, b_n\)\\
&&\qquad \le\ P\(\max_{1\le k\le n}\sum_{j=1}^{k} \f 1{j^{1/2}}
\(|X_j|I_{(|X_j|>\eta_j)}+E|X_j|I_{(|X_j|>\eta_j)}\)\ge C b_n/n^{1/2}\)\\
&&\qquad \le\ \f {2n^{1/2}}{C b_n}
\sum_{k=1}^{n}\f 1{k^{1/2}}E|X|I_{(|X|>\eta_k)}\\
&&\qquad \le\ \f {2n}{Cb_n}E|X|I_{(|X|>\eta_n)}+\f
{2C^{-1}}{l^{1/2}(\eta_n)}
\sum_{k=1}^{n}\f 1{k^{1/2}}E|X|I_{(\eta_k<|X|\le \eta_{n})}\\
&&\qquad \le\ AC^{-1}+\f {2C^{-1}}{l^{1/2}(\eta_n)}
\sum_{k=1}^{n}\eta^{-1/2}(\eta_k)\,E|X|^2I_{(\eta_k<|X|\le \eta_{k+1})}\\
&&\qquad \le \ AC^{-1}. \eestar This yields \be J_{2n}&:=&
b_n^{-1}\,\sup_{1/n\le t\le 1}
\Big|\sum_{j=1}^{[nt]}(X_j-X_j^*+EX_j^*)\Big|\Big/t^{1/2}=O_P(1).
\la {vic} \ee Combining the estimates (\ref {formu410}) and (\ref
{vic}), we obtain \bestar b_n^{-1} \sup_{0< t\le 1}
\Big|S_{[nt]}-\sum_{j=1}^{[nt]}\si_j^*Y_j\Big|\Big/t^{1/2}&\le &
J_{1n}+J_{2n}= O_P(1), \eestar which implies (\ref {formu4}). This also completes the
proof of Proposition \ref {pro1}.

\section{$L_p$-approximations of weighted self-normalized processes \la {sec5}}

Noting that $\limsup_{t\downarrow 0}|W(t)|/t^{1/2}=\infty\, a.s.$,
 it is impossible to extend Theorem \ref {th2} and its corollaries to
 the weight function
 $q(t)=t^{1/2}$. However,
  due to the finiteness of the integral
 $\int_0^1|W(t)|/t^{1/2}dt$, such a weight function is an immediate candidate for
weighted
 $L_1$-approximation.
 We still use $Q$ to denote   the class of those positive functions on $(0,1]$
for which $q(t)$  is nondecreasing near zero. The main purpose of
this section is to establish
 $L_p$-approximations of weighted self-normalized
 partial sum processes, and thus to extend Theorem 1.1 of Szyszkowicz (1993)
 that is based on assuming two moments when working with standardized partial sums.

\begin{th}\la {th3} Let $q\in Q$ and $0<p<\infty$.
Let $EX=0$ and $X$ be in the domain of attraction of the normal
law.
 Then the following statements hold true.
\begin{itemize}
\item[{\rm(a)}] On an appropriate probability space for the i.i.d. random variables
$X,X_1,X_2,...,$ we can construct a standard Wiener process
$\{W(t), 0\le t<\infty\}$
 such that, as $n\to\infty$,
\be 
\int_0^1\Big|V_n^{-1}S_{[nt]}-n^{-1/2}W(nt)\Big|^p\Big/q(t)dt
=o_P(1) \la {5th} \ee if and only if \be
\int_{0+}^1t^{p/2}/q(t)dt<\infty. \la {5th1} 
\ee

\item[{\rm(b)}] Let
$\{W(t), 0\le t<\infty\}$ be  a standard Wiener process. Then, as
$n\to\infty$, 
\be
\int_0^1\Big|V_n^{-1}S_{[nt]}\Big|^p\Big/q(t)dt\rightarrow_d
\int_0^1\Big|W(t)\Big|^p\Big/q(t)dt \la {5th2} 
\ee 
if and only if
$(\ref {5th1})$ holds.

\end{itemize}
\end{th}

It is of interest to call attention to the fact that (a) and (b)
are equivalent under the conditions of Theorem \ref {th3}. Thus,
unlike in the case of  Theorem \ref {th2} where, in addition to Corollary \ref {cor3}  we also have 
Corollary \ref {cor4} for a larger class of weight functions,in
$L_p$ we have (a) and (b) for the same class of weight functions.

In the proof of Theorem \ref {th3} we make use of the following
result, which is a consequence of Corollary 2.1 of Cs\"org\H{o},
Horv\'ath and Shao (1993).

\begin{lem} \la {lem8} Let $\{W(t), 0\le t<\infty\}$ be a standard Wiener process,
$0<p<\infty$, and assume that $q$ is a positive function on
$(0,1]$. Then $(\ref {5th1})$ holds if and only if \be
\int_{0+}^1|W(t)|^p/q(t)dt<\infty \qquad a.s. \la {sund} \ee
\end{lem}

We now turn to the proof of Theorem \ref {th3}.

\medskip
{\it Proof.} We first show  $(\mbox{a})$. If (\ref {5th}) is
satisfied, then
$$
\int_0^{1/n}\Big|n^{-1/2}W(nt)\Big|^p\Big/q(t)dt {=}_d
\int_0^{1/n}\Big|W(t)\Big|^p\Big/q(t)dt=o_P(1).
$$
Therefore, following the proof of Theorem 1.1 of Szyszkowicz
(1993), we have (\ref {sund}), and hence also (\ref {5th1}), via
Lemma \ref {lem8}.

We now prove that (\ref {5th1}) implies (\ref {5th}). Still  using
the notations in the proof of Theorem \ref {th1}, it can be easily
seen that for  $\delta\in (0,1),$
\newpage
\be
&&\int_0^1\Big|V_n^{-1}S_{[nt]}-n^{-1/2}W(nt)\Big|^p\Big/q(t)dt\no\\
&&\qquad \le \
\int_{\delta}^1\Big|V_n^{-1}S_{[nt]}-n^{-1/2}W(nt)\Big|^p\Big/q(t)dt\no\\
&&\qquad \quad +\ \(\f {b_n}{V_n}\)^p\int_{1/n}^{\delta}
\Big|b_n^{-1}\(S_{[nt]}-(S_{[nt]}^*-ES_{[nt]}^*)\)\Big|^p\Big/q(t)dt\no\\
&&\qquad \quad  +\ \(\f {b_n}{V_n}\)^p
\int_{1/n}^{\delta}\Big|b_n^{-1}\(S_{[nt]}^*-ES_{[nt]}^*\)-
b_n^{-1}\sum_{j=1}^{[nt]}\sigma_j^*Y_j\Big|^p\Big/q(t)dt\no\\
&&\qquad \quad +\ \(\f {b_n}{V_n}\)^p
\int_{1/n}^{\delta}\Big|b_n^{-1}
\sum_{j=1}^{[nt]}\sigma_j^*Y_j\Big|^p\Big/q(t)dt\ +\
\int_0^{\delta}\Big|n^{-1/2}W(nt)\Big|^p\Big/q(t)dt\no\\
&&\qquad :=\ II_1(n)+II_2(n)+II_3(n)+II_4(n)+II_5(n). \la {5th8}
\ee By Theorem \ref {th2}, we have \be II_1(n)&\le&\sup_{0\le t\le
1}\Big|V_n^{-1}S_{[nt]}-n^{-1/2}W(nt)\Big|^p
\int_{\delta}^11\Big/q(t)dt\no\\
&=&o_P(1), \quad \mbox{for any $\delta\in (0,1)$.} \la {5pro1} \ee
By (\ref {impor3}) and (\ref {vic}), it follows that \be
II_2(n)&\le &\(\f {b_n}{V_n}\)^p (J_{2n})^p
\int_{1/n}^{\delta}t^{p/2}\Big/q(t)dt=O_P(1)\int_{0}^{\delta}t^{p/2}\Big/q(t)dt.
\la {5vic} \ee Similarly, by (\ref{impor3}) and (\ref {formu410}),
we get \be II_3(n)&\le &\(\f {b_n}{V_n}\)^p (J_{1n})^p\,
\int_{1/n}^{\delta}t^{p/2}\Big/q(t)dt=O_P(1)\int_{0}^{\delta}t^{p/2}\Big/q(t)dt.
\la {5formu4} \ee Let $\delta$ be small enough so that $q$ is
already nondecreasing on $(0,\delta)$. Similarly  to the proof of
(\ref {addtion}), we have 
\bestar
\int_{1/n}^{\delta}\Big|b_n^{-1}\sum_{j=1}^{[nt]}\sigma_j^*Y_j\Big|\Big/q(t)dt
&{=}_d & \int_{1/n}^{\delta}\Big|W\Big(\f
1n\sum_{j=1}^{[nt]}
\f {\sigma_j^{*2}}{l(\eta_n)}\Big)\Big|\Big/q(t)dt\no\\
&\le & \int_{1/n}^{\delta}\Big|W\Big(\f 1n\sum_{j=1}^{[nt]} \f
{\sigma_j^{*2}}{l(\eta_n)}\Big)\Big|\Big/q\Big(\f
1n\sum_{j=1}^{[nt]} \f {\sigma_j^{*2}}{l(\eta_n)}\Big)dt. \eestar
This, together with (\ref {impor3}), implies that \be
II_4(n)&=&O_P(1)\int_{0}^{\delta}\Big|W\Big(\f 1n\sum_{j=1}^{[nt]}
\f {\sigma_j^{*2}}{l(\eta_n)}\Big)\Big|\Big/q\Big(\f
1n\sum_{j=1}^{[nt]} \f {\sigma_j^{*2}}{l(\eta_n)}\Big)dt. 
\la{5addtion} 
\ee 
We also have for each $n$ 
\be 
II_5(n)&{=}_d &\int_0^{\delta}\Big|W(t)\Big|^p\Big/q(t)dt. \la
{5formu5} \ee Taking $n\to\infty$ and then $\delta$ arbitrarily
small, on using (\ref {5th1}) via Lemma \ref {lem8}, by (\ref
{5th8})-(\ref {5formu5}) we arrive at  (\ref {5th}).

We next show  (b). If (\ref {5th2}) holds true, then clearly we
have (\ref {sund}) and hence also (\ref {5th1}) by Lemma \ref
{lem8}. Conversely, we now assume  (\ref {5th1}). Since $q$ is a
positive function on $(0,1]$, it follows from  Theorem \ref {th2}
that, as $n\to\infty$,
$$
\int_{\delta}^1\Big|V_n^{-1}S_{[nt]}\Big|^p\Big/q(t)dt\rightarrow_d
\int_{\delta}^1\Big|W(t)\Big|^p\Big/q(t)dt
$$
for all $\delta\in (0, 1)$. On using (\ref {5th1}), by Lemma \ref
{lem8} we conclude that
$$
\lim_{\delta\downarrow
0}\int_0^{\delta}\Big|W(t)\Big|^p\Big/q(t)dt=0\qquad a.s.
$$
From here on, mutatis mutandis, the proof follows similarly to that of proving (iii)
of Theorem 1.1 of Cs\"org\H{o}, Horv\'ath and Shao (1993).

This  completes the proof of Theorem \ref {th3}.

\vspace{.5cm}
\begin{center}
{\large{\bf REFERENCES}}
\end{center}
\baselineskip=0.5 true cm

\newcommand{
\itemr}[1] {\par\indent\hbox to -29.0pt{#1\hfil} \hangindent=20.0pt
\hangafter=1\hskip 10pt}

\itemr{} Anderson, T. W. and Darling, D. A.~(1952). Asymptotic theory of certain
``goodness of fit" criteria based on stochastic processes, {\it
Ann. Math. Statist.} {\bf 23}, 193-212.

\smallskip
\itemr{} Breiman, L.~(1968).  {\it Probability}, Addison-Wesley, Reading, MA.

\smallskip\itemr{} Chibisov, D.~(1964). Some theorems on the limiting behaviour of
empirical distribution functions, {\it  Selected Transl. Math.
Statist. Probab.} {\bf 6}, 147-156.

\smallskip\itemr{}
 Chistyakov, G. P. and G\"otze, F. (2004).
  Limit distributions of Studentized means. {\it Ann. Probab.} {\bf 32}, no. 1A,
28--77.

\smallskip\itemr{} Cs\"org\H{o}, M.~(2002).  A glimpse of the impact of P\'al Erd\H{o}s on probability and statistics,  {\it The Canadian Journal of Statistics} {\bf 30}, 493--556.

\smallskip\itemr{} Cs\"org\H{o}, M., Cs\"org\H{o}, S., Horv\'ath, L. and Mason, D.M.~(1986).
Weighted empirical and quantile processes, {\it Ann. Probab.} {\bf
14}, 31-85.

\smallskip\itemr{} Cs\"org\H{o}, M. and  Horv\'ath, L.~(1986). Approximations of weighted
empirical and quantile processes, {\it Statist.\ Probab.\ Lett.} {\bf 4}, 151--154.

\smallskip\itemr{} Cs\"org\H{o}, M. and  Horv\'ath, L.~(1988). Nonparametric methods
for changepoint problems, In {\it Handbook of Statistics}, {\bf
7}, Elsevier Science Publisher B.V.,  403-425, North-Holland,
Amsterdam.

\smallskip\itemr{} Cs\"org\H{o}, M. and  Horv\'ath, L.~(1993).
{\it Weighted Approximations in Probability and Statistics}, Wiley
Series in Probability and Mathematical Statistics: Probability and
Mathematical Statistics, Wiley, Chichester.

\smallskip\itemr{} Cs\"org\H{o}, M. and  Horv\'ath, L.~(1997). {\it Limit Theorems
in Change-Point Analysis}, Wiley Series in Probability and
Mathematical Statistics: Probability and Mathematical Statistics,
Wiley, Chichester.

\smallskip\itemr{} Cs\"org\H{o}, M. and  Horv\'ath, L.\ and Shao, Q.-M.~(1993). Convergence of
integrals of uniform empirical and quantile processes, {\it Stochastic Process.\
Appl.} {\bf 45}, 283--294.

\smallskip\itemr{} Cs\"org\H{o}, M., Norvai\v{s}a, R.\ and Szyszkowicz, B.~(1999).
Convergence of weighted partial sums when the limiting
distribution is not necessarily Radon, {\it Stochastic Process.
Appl.} {\bf 81}, 81-101.

\smallskip\itemr{} Cs\"org\H{o}, M. and R\'ev\'esz, P.~(1981).  {\it Strong Approximations in Probability and Statistics}, Academic Press, New York.

\smallskip\itemr{} Cs\"org\H{o}, M., Shao, Q.-M.\  and  Szyszkowicz, B.\ (1991). A note on
local and global functions of a Wiener process and some R\'enyi-type statistics,
{\it Studia Sci.\ Math.\ Hungar.} {\bf 26}, 239--259.

\smallskip\itemr{} Cs\"org\H{o}, M., Szyszkowicz, B. and Wang, Q.~(2003).
Donsker's theorem for self-normalized partial sums processes. {\it
Ann. Probab.} {\bf 31}, 1228--1240.

\smallskip\itemr{}
Cs\"org\H{o}, M., Szyszkowicz, B. and  Wang, Q. (2004). On
weighted approximations and strong limit theorems for
self-normalized partial sums processes. In {\it Asymptotic Methods
in Stochastics,}  489--521, Fields Inst. Commun. {\bf 44}, Amer.
Math. Soc., Providence, RI.

\smallskip\itemr{} Gine, E., G\"otze, F. and Mason, D.M.~(1997).
When is the Student $t$-statistic asymptotically standard normal?
{\it Ann. Probab.} {\bf 25}, 1514-1531.

\smallskip\itemr{}Griffin, P.S.\ (2002).  Tightness of the Student t--statistic, {\it
Electron.\ Comm.
 Probab.} {\bf 2}, 181-190.

\smallskip\itemr{}
 Koml\'os, J., Major, P.\ and Tusn\'ady, G.~(1975).  An approximation of partial sums of independent RV's, and the sample DF, I., {\it Z.~Wahrscheinlichkeitstheorie verw.\ Gebiete} {\bf 32}, 111--131.

\smallskip\itemr{}
 Koml\'os, J., Major, P.\ and Tusn\'ady, G.
 (1976). An approximation of partial sums of independent RV's, and the sample DF. II.,
 {\it Z. Wahrscheinlichkeitstheorie verw. Gebiete} {\bf 34},  33--58.

\smallskip\itemr{} Logan, B.F., Mallows, C.L., Rice, S.O.\ and Shepp, L.A.~(1973), Limit
distributions of self-normalized sums.  {\it Ann.\ Probab.} {\bf 1}, 788--809.
\smallskip\itemr{} O'Reilly, P.~(1974). On the weak convergence of empirical processes in
sup-norm metric,
{\it Ann. Probab.} {\bf 2}, 642-651.
\smallskip\itemr{}
 Major, P.(1976). The approximation of partial sums of independent RV's. {\it Z.
Wahrscheinlichkeitstheorie und Verw. Gebiete} {\bf 35}, 213--220.

\smallskip\itemr{} Pyke, R. and Shorack, G.R.~(1968). Weak convergence of two-sample
 empirical processes and a new approach to Chernoff-Savage theorems, {\it Ann. Math.
Statist.}
 {\bf 39}, 755-771.

\smallskip\itemr{} R\'enyi, A.~(1953). On the theory of order statistics, {\it Acta. Math.
Acad. Sci.
Hungar.} {\bf 4}, 191-231.

\smallskip\itemr{} Sakhanenko, A. I.~(1980). On unimprovable estimates of the rate of convergence
in invariance principle, In {\it Colloquia Math. Soc. J\'anos
Bolyai {\bf 32}, Nonparametric Statistical Inference,} Budapest,
Hungary, 779-783, North-Holland, Amsterdam.

\smallskip\itemr{} Sakhanenko, A. I.~(1984). On  estimates of the rate of convergence
in the invariance principle, In {\it Advances in Probability
Theory:
 Limit Theorems and Related Problems,} (A. A. Borovkov, Ed.) 124-135. Springer, New
York.

\smallskip\itemr{} Sakhanenko, A. I.~(1985). Convergence rate in the invariance principle for
 non-identically distributed variables with exponential moments,
  In {\it Advances in Probability Theory:
 Limit Theorems for Sums of  Random Variables,} (A. A. Borovkov, Ed.) 2-73.
  Springer, New York.

\smallskip\itemr{} Shorack, G. R. and Wellner, J. A.~(1986).
{\it Empirical Processes with Applications to Statistics}, Wiley
Series in Probability and Mathematical Statistics: Probability and
Mathematical Statistics, Wiley, New York.

\smallskip\itemr{} ``Student" (1908).  The probable error of the mean, {\it Biometrika} {\bf
6}, 1-25.
\smallskip\itemr{} Szyszkowicz, B.~(1991). Weighted stochastic processes under contiguous
alternatives,  {\it C.R. Math. Rep. Acad. Sci. Canada } {\bf 13},
211-216.

\smallskip\itemr{} Szyszkowicz, B.~(1992). Weighted asymptotics of partial sum processes in
$D[0,\infty)$,
  {\it C.R. Math. Rep. Acad. Sci. Canada } {\bf 14}, 273-278.

\smallskip\itemr{} Szyszkowicz, B.~(1993). $L_p$-approximations of weighted
partial sum processes, {\it Stochastic Process. Appl.} {\bf 45},
295-308.

\smallskip\itemr{} Szyszkowicz, B.~(1996). Weighted approximations of partial sum processes
in $D[0,\infty)$. I, {\it Studia Sci. Math. Hungar.} {\bf 31},
323-353.

\smallskip\itemr{} Szyszkowicz, B.~(1997). Weighted approximations of partial sum processes
in $D[0,\infty)$. II, {\it Studia Sci. Math. Hungar.} {\bf 33},
305-320.

\vskip1cm

\noindent
Mikl\'os Cs\"org\H{o} \\
School of Mathematics and Statistics\\
Carleton University\\
1125 Colonel By Drive\\
Ottawa, ON ~Canada K1S 5B6\\
{\tt mcsorgo@math.carleton.ca}
\\[0.5cm]
Barbara Szyszkowicz \\
School of Mathematics and Statistics\\
Carleton University\\
1125 Colonel By Drive\\
Ottawa, ON ~Canada K1S 5B6\\
{\tt bszyszko@math.carleton.ca}
\\[0.5cm]
Qiying Wang \\
School of Mathematics and Statistics\\
University of Sydney\\
NSW 2006, Australia\\
{\tt qiying@maths.usyd.edu.au}

\end{document}